\documentclass{amsart}
\usepackage{graphicx}
\usepackage{amssymb}
\usepackage{epstopdf}
\usepackage{nicefrac}
\usepackage{enumerate}
\usepackage{hyperref}
\usepackage{float}
\usepackage{color}
\usepackage{pst-all}
\usepackage{tabularx}

\newtheorem{thm}{THEOREM}[section]

\newtheorem{cor}[thm]{COROLLARY}

\newtheorem{defn}[thm]{DEFINITION}
\newtheorem{ex}[thm]{EXAMPLE}

\newtheorem{prop}[thm]{PROPOSITION}
\newtheorem{quest}[thm]{QUESTION}

\newcommand{\ds}{\displaystyle}

\newcommand{\cG}{{\mathcal G}}
\newcommand{\cH}{{\mathcal H}}

\newcommand{\cO}{{\mathcal O}}

\newcommand{\CO}{{\rm CO}}

\newcommand{\cU}{{\mathcal U}}

\newcommand{\diam}{{\rm diam}} 
 
\newcommand{\dX}{d_{\fX}} 
\newcommand{\e}{{\varepsilon}} 

\newcommand{\fD}{{\mathfrak{D}}}
\newcommand{\fH}{{\mathfrak{H}}}
\newcommand{\fG}{{\mathfrak{G}}}

\newcommand{\fN}{{\mathfrak{N}}}

\newcommand{\fX}{{\mathfrak{X}}}

\newcommand{\fU}{{\mathfrak{U}}}
\newcommand{\G}{\Gamma}

\newcommand{\Homeo}{{\rm Homeo}}

\newcommand{\mZ}{{\mathbb Z}}

\newcommand{\vp}{{\varphi}}
\newcommand{\whe}{\widehat{e}}
\newcommand{\whg}{\widehat{g}}

\newcommand{\whC}{\widehat{C}}

\newcommand{\whU}{\widehat{U}}
\newcommand{\whmZ}{\widehat{\mZ}}
\newcommand{\whPhi}{\widehat{\Phi}}

\newcommand{\whtau}{{\widehat{\tau}}}

\newcommand{\whGamma}{\widehat{\Gamma}}
\newcommand{\whpi}{\widehat{\pi}}

\newcommand{\oa}{\overline{a}}
\newcommand{\ob}{\overline{b}}
\newcommand{\oc}{\overline{c}}

\newcommand{\ox}{\overline{x}}
\newcommand{\oy}{\overline{y}}
\newcommand{\oz}{\overline{z}}

\newcommand{\wha}{\widehat{a}}
\newcommand{\whb}{\widehat{b}}
\newcommand{\whc}{\widehat{c}}

\newcommand{\lcm}{{\rm lcm}}

\newcommand\mor{\mathrel{\stackrel{\makebox[0pt]{\mbox{\normalfont\tiny a}}}{\sim}}}

\begin{document}

\title{Prime spectrum and dynamics for nilpotent Cantor actions}

\author{Steven Hurder}
\address{Steven Hurder, Department of Mathematics, University of Illinois at Chicago, 322 SEO (m/c 249), 851 S. Morgan Street, Chicago, IL 60607-7045}
\email{hurder@uic.edu}

\author{Olga Lukina}
\address{Olga Lukina, Mathematical Institute, Leiden University, PO Box 9512, 2300 RA Leiden, The Netherlands}
\email{o.lukina@math.leidenuniv.nl}

\thanks{Version date: May 1, 2023}

\thanks{2020 {\it Mathematics Subject Classification}. Primary: 20E18, 37B05, 37B45; Secondary: 57S10}

\thanks{Keywords: odometers, Cantor actions, profinite groups, Steinitz numbers, Heisenberg group}

  \begin{abstract}
  A minimal equicontinuous action by homeomorphisms of a discrete group $\Gamma$ on a Cantor set $X$   is locally quasi-analytic, if each homeomorphism  has a unique extension  from small open sets to open sets of uniform diameter on $X$. A minimal action is stable, if the actions of   $\Gamma$  and of the closure of $\Gamma$ in the group of homeomorphisms of $X$, are both locally quasi-analytic.  
  
When $\Gamma$ is virtually nilpotent, we say that $\Phi \colon \Gamma \times \mathfrak{X} \to \mathfrak{X}$ is a nilpotent Cantor action.    We show that a nilpotent Cantor action with finite prime spectrum must be stable.  We also  prove   there exist uncountably many distinct   Cantor actions of the Heisenberg group, necessarily with infinite prime spectrum, which are not stable. 
  \end{abstract}

\maketitle
    
\section{Introduction}\label{sec-intro}

 A minimal equicontinuous action $\Phi \colon \G \times \fX \to \fX$ of a countable group $\G$ on a Cantor space $\fX$ is called a \emph{generalized odometer} \cite{CortezPetite2008,GPS2019}. When $\G = \mZ$, this is just the abstract form of a traditional odometer action of the integers. For $\G = \mZ^n$ with $n \geq 2$, one obtains a more complex class of actions, whose classification becomes increasingly intractable as $n$ increases \cite{Thomas2003}, even while the dynamical properties of minimal equicontinuous Cantor actions by $\mZ^n$ are well-behaved. For   $\G$ in general, we simply refer to these as \emph{Cantor actions}, which will always be assumed \emph{minimal} and \emph{equicontinuous}.

It is a classical result that a $\mZ$-odometer is classified by its Steinitz order, which is calculated using a representation of the action as an inverse limit of actions on finite cyclic groups. One can also associate to a   Cantor action by $\mZ^n$ its Steinitz order, and also a collection of types, called its \emph{typeset}, which consists of equivalence classes of Steinitz orders of individual elements of $\mZ^n$. As discussed by Thomas \cite[Section~4]{Thomas2008}, the additional data of the typeset is still not sufficient to reduce the classification problem for Cantor actions by $\mZ^n$ to a standard Borel equivalence relation.

In the authors' work \cite{HL2023a}, we associate the \emph{type} and \emph{typeset} invariants  to a Cantor action $(\fX, \G , \Phi)$ for arbitrary countable group $\G$. The type $\tau[\fX, \G , \Phi]$ is the asymptotic equivalence class of the Steinitz order $\xi(\fX, \G , \Phi)$ of a presentation of the action as an inverse limit of actions of 
$\G$ on finite sets.  

Associated to the type $\tau[\fX, \G , \Phi]$ is an even more basic invariant, the \emph{prime spectrum} $\pi[\fX, \G , \Phi]$, which consists of the set of primes which appear   in a Steinitz order $\xi(\fX, \G , \Phi)$ representing the type $\tau[\fX, \G , \Phi]$, see Section \ref{def-primespectrum}. The prime spectrum decomposes into two parts, 
$$\pi[\fX, \G , \Phi] = \pi_{\infty}[\fX, \G , \Phi] \cup \pi_f[\fX, \G , \Phi],$$ 
where  the \emph{infinite prime spectrum} $\pi_{\infty}[\fX, \G , \Phi]$ consists of the primes which occur with infinite multiplicity in $\xi(\fX, \G , \Phi)$, and the \emph{finite prime spectrum} $\pi_f[\fX, \G , \Phi]$ which consists of the primes that occur with finite multiplicity.  
The prime spectrum and the finite prime spectrum are only well-defined  modulo finite subsets of $\pi_f[\fX, \G , \Phi]$.  
\begin{defn}
A Cantor action $(\fX, \G , \Phi)$ has \emph{finite spectrum} if the prime spectrum $\pi[\fX, \G , \Phi]$ is a finite set, and is said to have  \emph{infinite   spectrum} otherwise. 
\end{defn}
The  classification of  Cantor actions for $\G$ is, in general,  intractable and one seeks invariants for Cantor actions which at least distinguish between particular classes of actions.    The authors' works \cite{HL2019,HL2020,HL2021,HL2022} study   dynamical properties which yield invariants of Cantor actions.  In particular,   one of the most basic invariants, is the property that the action is \emph{stable} or \emph{wild}.   The purpose of this note is to  give a relation between the prime spectrum of a Cantor action and the     wild property. 

As explained in detail in Section~\ref{subsec-lqa} below, the property that the action $(\fX,\G,\Phi)$ is stable is a property of the action of the completion $\fG(\Phi) = \overline{\Phi(\G)} \subset {\rm Homeo}(\fX)$, which is a profinite group naturally acting on $\fX$. The property that the action $(\fX,\G,\Phi)$ is locally quasi-analytic is defined in Definition \ref{def-LQA}, and  $(\fX,\G,\Phi)$ is \emph{stable} if the action of $\fG(\Phi)$ on $\fX$ is also locally quasi-analytic. If $(\fX,\G,\Phi)$ is stable, then $(\fX,\G,\Phi)$ is locally quasi-analytic. The converse need not hold even for actions of nilpotent groups, as we see later.

 A Cantor action $(\fX, \G , \Phi)$  is said to be \emph{nilpotent} if  $\G$ contains a finitely generated nilpotent subgroup with finite index. This class of group actions is particularly interesting, as it has  the natural next level of complexity after the abelian Cantor actions. 
  We show the following three results  for   nilpotent Cantor actions.

    \begin{thm}\label{thm-main1}
Let $(\fX,\G,\Phi)$ be a nilpotent Cantor action. 
If the prime spectrum $\pi[\fX,\G,\Phi]$ is   finite,  then the action is stable.
\end{thm}

Theorem \ref{thm-main1} does not have a converse. We   show that every collection of primes, finite or infinite, can be realized as the prime spectrum of a stable nilpotent Cantor action.

\begin{thm}\label{thm-main3}
Let $\pi_f$ and $\pi_\infty$ be two distinct sets of primes, where $\pi_f$ is a finite set, and $\pi_\infty$ is a non-empty finite or infinite set. Then there exists a stable nilpotent Cantor action $(\fX,\G,\Phi)$ such that $ \pi_{\infty}[\fX, \G , \Phi] = \pi_\infty$ and $\pi_f[\fX, \G , \Phi] = \pi_f$.
\end{thm}

Let $(\fX,\G,\Phi)$ be an Cantor action. If the action is effective,   then it is free, and   the action of the closure $\fG(\Phi)$ is also free, which implies that the action is stable. An effective nilpotent Cantor action need not be free, and may even have elements which fix every point in a clopen subset of the Cantor set $\fX$. The authors showed in their work \cite{HL2020} that nilpotent Cantor actions are locally quasi-analytic, which means that such subsets of fixed points cannot be arbitrarily small, as their diameter has lower bound which is uniform over the Cantor set $\fX$.  It is then  surprising to discover that if  one allows $\fG(\Phi)$ to have infinite prime spectrum then one can construct wild nilpotent actions, for which the action of the closure $\fG(\Phi)$ is not locally quasi-analytic, as shown in Theorem \ref{thm-main2}. 
In addition, Theorem \ref{thm-main2} is a realization result, which shows that every infinite set of primes can be realized as the prime spectrum of a wild nilpotent Cantor action.

    \begin{thm}\label{thm-main2}
Given any two distinct sets $\pi_f$ and $\pi_\infty$ of primes, where $\pi_f$ is infinite and $\pi_\infty$ is any (possibly empty) set, there is a minimal equicontinuous action $(\fX,\G,\Phi)$ of the Heisenberg group, such that $\pi_f[\fX,\G,\Phi] = \pi_f$ and $\pi_\infty[\fX,\G,\Phi] = \pi_\infty$.

 Moreover, there exists an uncountable number of nilpotent Cantor actions $(\fX,\G,\Phi)$ of the Heisenberg group $\G$ with infinite prime spectra such that:
  \begin{enumerate}
  \item Each $(\fX,\G,\Phi)$ is topologically free,
  \item Each $(\fX,\G,\Phi)$ is wild,
  \item The prime spectra of such actions are pairwise distinct.
  \end{enumerate}
    \end{thm} 
    
  The notion of return equivalence for Cantor actions, and its relationship with conjugacy of action is explained in Section \ref{subsec-equivalence}. The result of Corollary \ref{cor-main4} below follows from the result that the prime spectrum of the action is an invariant of its return equivalence class, see Theorem \ref{thm-welldefined}.
 
\begin{cor}\label{cor-main4}
There exists an uncountable number of nilpotent Cantor actions $(\fX,\G,\Phi)$ of the Heisenberg group $\G$ which are not return equivalent, and therefore not conjugate.
\end{cor}

The conclusion of Theorem~\ref{thm-main2} is used in   \cite{HL2023c} for the calculation of the mapping class groups  of solenoidal manifolds whose base is a nil-manifold.

We note that for more general groups $\G$, an analog of Theorem~\ref{thm-main1} need not hold. For example, a weakly branch group, as studied in \cite{Bartholdi2003,BartholdiGrigorchuk2002,BGS2012,Nekrashevych2005}, acts    on the boundary of a $d$-regular rooted tree, and so   has finite prime spectrum $\{d\}$, but the dynamic of the actions on the Cantor   boundary are wild.  
\begin{quest}\label{quest-types}
Let $(\fX,\G,\Phi)$ be a   Cantor action. For which classes of groups $\G$ does the finiteness of the prime spectrum of the action imply that the action is stable?
\end{quest}

The paper is organized as follows. In Section \ref{subsec-basics} we recall basic properties of minimal equicontinuous group actions on Cantor sets. In particular, the definition of the prime spectrum of a minimal equicontinuous action is given in Definition \ref{def-primespectrum}. We prove Theorem \ref{thm-main1} in Section \ref{sec-nilpotent}, and give basic examples of nilpotent Cantor actions in Section \ref{sec-examples}. In Section \ref{sec-prescribed} we construct examples of stable and wild actions of the Heisenberg group with prescribed prime spectrum, proving Theorems \ref{thm-main3} and \ref{thm-main2}, from which we deduce  Corollary \ref{cor-main4}.

 \section{Cantor actions}\label{sec-basics}

We recall some of the basic 
 properties of     Cantor actions, as required for the proofs of Theorems~\ref{thm-main1} and \ref{thm-main2}.
More complete discussions of the properties of equicontinuous Cantor actions are given in     the text by Auslander \cite{Auslander1988}, the papers by Cortez and Petite  \cite{CortezPetite2008}, Cortez and Medynets  \cite{CortezMedynets2016},    and   the authors' works, in particular   \cite{DHL2016a,DHL2016c} and   \cite[Section~3]{HL2021}.

\subsection{Basic concepts}\label{subsec-basics}

Let  $(\fX,\G,\Phi)$   denote an action  $\Phi \colon \G \times \fX \to \fX$. We   write $g\cdot x$ for $\Phi(g)(x)$ when appropriate.
The orbit of  $x \in \fX$ is the subset $\cO(x) = \{g \cdot x \mid g \in \G\}$. 
The action is \emph{minimal} if  for all $x \in \fX$, its   orbit $\cO(x)$ is dense in $\fX$.

 An action  $(\fX,\G,\Phi)$ is \emph{equicontinuous} with respect to a metric $\dX$ on $\fX$, if for all $\e >0$ there exists $\delta > 0$, such that for all $x , y \in \fX$ and $g \in \G$ we have  that 
 $\ds  \dX(x,y) < \delta$ implies   $\dX(g \cdot x, g \cdot y) < \e$.
The property of being equicontinuous    is independent of the choice of the metric   on $\fX$ which is   compatible with the topology of $\fX$.

Now assume that $\fX$ is a Cantor space. 
Let $\CO(\fX)$ denote the collection  of all clopen (closed and open) subsets of  $\fX$, which forms a basis for the topology of $\fX$. 
For $\phi \in \Homeo(\fX)$ and    $U \in \CO(\fX)$, the image $\phi(U) \in \CO(\fX)$.  
The following   result is folklore, and a proof is given in \cite[Proposition~3.1]{HL2020}.
 \begin{prop}\label{prop-CO}
 For $\fX$ a Cantor space, a minimal   action   $\Phi \colon \G \times \fX \to \fX$  is  equicontinuous  if and only if  the $\G$-orbit of every $U \in \CO(\fX)$ is finite for the induced action $\Phi_* \colon \G \times \CO(\fX) \to \CO(\fX)$.
\end{prop}
 
 \begin{defn} 
We say that $U \subset \fX$  is \emph{adapted} to the action $(\fX,\G,\Phi)$ if $U$ is a   \emph{non-empty clopen} subset, and for any $g \in \G$, 
if $\Phi(g)(U) \cap U \ne \emptyset$ implies that  $\Phi(g)(U) = U$.   
\end{defn}

The proof of   \cite[Proposition~3.1]{HL2020} shows that given  $x \in \fX$ and clopen set $x \in W$, there is an adapted clopen set $U$ with $x \in U \subset W$. 

For an adapted set $U$,   the set of ``return times'' to $U$, 
 \begin{equation}\label{eq-adapted}
\G_U = \left\{g \in \G \mid g \cdot U  \cap U \ne \emptyset  \right\}  
\end{equation}
is a subgroup of   $\G$, called the \emph{stabilizer} of $U$.      
  Then for $g, g' \in \G$ with $g \cdot U \cap g' \cdot U \ne \emptyset$ we have $g^{-1} \, g' \cdot U = U$, hence $g^{-1} \, g' \in \G_U$. Thus,  the  translates $\{ g \cdot U \mid g \in \G\}$ form a finite clopen partition of $\fX$, and are in 1-1 correspondence with the quotient space $X_U = \G/\G_U$. Then $\G$ acts by permutations of the finite set $X_U$ and so the stabilizer group $\G_U \subset G$ has finite index.  Note that this implies that if $V \subset U$ is a proper inclusion of adapted sets, then the inclusion $\G_V \subset \G_U$ is also proper.

\begin{defn}\label{def-adaptednbhds}
Let  $(\fX,\G,\Phi)$   be a   Cantor    action.
A properly descending chain of clopen sets $\cU = \{U_{\ell} \subset \fX  \mid \ell > 0\}$ is said to be an \emph{adapted neighborhood basis} at $x \in \fX$ for the action $\Phi$,   if
    $x \in U_{\ell +1} \subset U_{\ell}$  is a proper inclusion for all $ \ell > 0$, with     $\cap_{\ell > 0}  \ U_{\ell} = \{x\}$, and  each $U_{\ell}$ is adapted to the action $\Phi$.
\end{defn}
Given $x \in \fX$ and   $\e > 0$, Proposition~\ref{prop-CO} implies there exists an adapted clopen set $U \in \CO(\fX)$ with $x \in U$ and $\diam(U) < \e$.  Thus, one can choose a descending chain $\cU$ of adapted sets in $\CO(\fX)$ whose intersection is $x$, from which the following result follows:

\begin{prop}\label{prop-adpatedchain}
Let  $(\fX,\G,\Phi)$   be a  Cantor    action. Given $x \in \fX$, there exists an adapted neighborhood basis $\cU$ at $x$ for the action $\Phi$.
 \end{prop}
 Combining the above remarks, we have:
\begin{cor}\label{cor-Uchain}  
Let  $(\fX,\G,\Phi)$   be a  Cantor    action, and   $\cU$ be an adapted neighborhood basis. Set  $\G_{\ell} = \G_{U_{\ell}}$,  with $\G_0 = \G$, then there is     a nested chain of finite index subgroups,  
$\ds  \cG_{\cU} = \{\G_0 \supset \G_1 \supset \cdots \}$.
 \end{cor}

\subsection{Profinite completion}\label{subsec-profinitecomp}

Let  $\Phi(\G) \subset \Homeo(\fX)$ denote the image subgroup for an action $(\fX, \G, \Phi)$.
When the action is equicontinuous,     the  closure $\overline{\Phi(\G)} \subset \Homeo(\fX)$ in the \emph{uniform topology of maps} is a separable profinite group.  We adopt the notation $\fG(\Phi) \equiv \overline{\Phi(\G)}$.    

Let $\whPhi \colon \fG(\Phi) \times \fX \to \fX$ denote the induced   action of $\fG(\Phi)$ on $\fX$, which is transitive as the action   $(\fX, \G, \Phi)$ is   minimal. For $\whg \in \fG(\Phi)$, we   write its action on $\fX$ by $\whg \, x = \whPhi(\whg)(x)$.
Given $x \in \fX$,   introduce the isotropy group,  
\begin{align}\label{iso-defn2}
 \fD(\Phi, x) = \{ \whg  \in \fG(\Phi) \mid \whg \, x = x\} \subset \Homeo(\fX) \ ,
\end{align}
which is a closed subgroup of $\fG(\Phi)$, and    thus is either finite, or is an infinite profinite group.    
 As the action $\whPhi \colon \fG(\Phi) \times \fX \to \fX$ is transitive,  the conjugacy class of $\fD(\Phi,x)$ in $\fG(\Phi)$ is independent of the choice of   $x$, and by abuse of notation we omit the subscript $x$.     The group $\fD(\Phi)$ is called the \emph{discriminant} of the action $(\fX,\G,\Phi)$ in  \cite{DHL2016c,HL2019,HL2021}, and is called a \emph{parabolic} subgroup  (of the profinite completion of a countable group)  in the works by Bartholdi and Grigorchuk \cite{BartholdiGrigorchuk2000,BartholdiGrigorchuk2002}.

 \subsection{Algebraic Cantor  actions}\label{subsec-gchains}
 
 Let $\cG = \{\G = \G_0 \supset \G_1 \supset \G_2 \supset \cdots\}$ be a descending chain of finite index subgroups.
 Let $X_{\ell} = \G/\G_{\ell}$ and note that  $\G$ acts transitively on the left on the finite set  $X_{\ell}$.    
The inclusion $\G_{\ell +1} \subset \G_{\ell}$ induces a natural $\G$-invariant quotient map $p_{\ell +1} \colon X_{\ell +1} \to X_{\ell}$.
 Introduce the inverse limit 
  \begin{eqnarray} 
X_{\infty} & \equiv &  \lim_{\longleftarrow} ~ \{ p_{\ell +1} \colon X_{\ell +1} \to X_{\ell}  \mid \ell \geq 0 \} \label{eq-invlimspace}\\
& = &  \{(x_0, x_1, \ldots ) \in X_{\infty}  \mid p_{\ell +1 }(x_{\ell + 1}) =  x_{\ell} ~ {\rm for ~ all} ~ \ell \geq 0 ~\} ~ \subset \prod_{\ell \geq 0} ~ X_{\ell} \  .  \nonumber
\end{eqnarray}
Then $X_{\infty}$  is a Cantor space with the Tychonoff topology, where the actions of $\G$ on the factors $X_{\ell}$ induce    a minimal  equicontinuous action  denoted by  $\Phi  \colon \G \times X_{\infty} \to X_{\infty}$. 
  There is a natural basepoint $x_{\infty} \in X_{\infty}$ given by the cosets of the identity element $e \in \G$, so $x_{\infty} = (e \G_{\ell})$. An adapted neighborhood basis   of $x_{\infty}$ is given by the clopen sets 
\begin{equation}\label{eq-openbasis}
U_{\ell} = \left\{ x = (x_{i}) \in X_{\infty}   \mid  x_i = e \G_i \in X_i~, ~ 0 \leq i \leq \ell ~  \right\} \subset X_{\infty} \ .
\end{equation}
There is a tautological identity $\G_{\ell} = \G_{U_{\ell}}$.

 For each $\ell \geq 0$, we have the ``partition coding map'' $\Theta_{\ell} \colon \fX \to X_{\ell}$ which is $\G$-equivariant.  The maps $\{\Theta_{\ell}\}$ are compatible with the   map on quotients in \eqref{eq-invlimspace}, and so they induce a  limit map $\Theta_x \colon \fX \to X_{\infty}$. The fact that the diameters of the clopen sets $\{U_{\ell}\}$ tend to zero, implies that $\Theta_x$ is a homeomorphism.  Moreover, $\Theta_x(x) =  x_{\infty} \in X_{\infty}$ where   $\{x\} = \cap_{\ell > 0}  \ U_{\ell}$.  
 The following is folklore:
 \begin{thm}\cite[Appendix~A]{DHL2016a} \label{thm-algmodel}
 The map $\Theta_x \colon  \fX \to X_{\infty}$ induces an isomorphism of the Cantor actions $(\fX,\G,\Phi)$ and $(X_{\infty}, \G, \Phi_x)$.
 \end{thm}
 The   action $(X_{\infty}, \G, \Phi_x)$   is called the \emph{odometer model} centered at $x$ for the action $(\fX,\G,\Phi)$.
 The dependence of the model   on the choices of a base point $x \in \fX$ and adapted neighborhood basis $\cU$ is discussed in detail in the   works \cite{DHL2016a,FO2002,HL2019,HL2021}.   Again, we abuse notation in the following and   omit the subscript ``$x$''.

 Next, we develop the algebraic model for the profinite action  $\whPhi \colon \fG(\Phi) \times \fX \to \fX$ of the completion  $\fG(\Phi) \equiv \overline{\Phi(\G)} \subset \Homeo(\fX)$. 
Fix a choice of group chain $\{\G_{\ell} \mid \ell \geq 0\}$ as above, which provides an algebraic model for the action $(\fX,\G,\Phi)$.

 For each $\ell \geq 1$, let $C_{\ell} \subset \G_{\ell}$ denote the \emph{core} of   $\G_{\ell}$, i.e. the largest normal subgroup of   $\G_{\ell}$ in $\G$. So 
\begin{equation}\label{eq-core}
C_{\ell} ~ = {\rm Core}(\G_{\ell}) ~ = ~ \bigcap_{g \in \G} ~ g \ \G_{\ell} \ g^{-1} ~ \subset \G_{\ell} ~ .
\end{equation}
As   $\G_{\ell}$ has finite index in $\G$, the same holds for $C_{\ell}$. Observe that for all $\ell \geq 0$,   we have $C_{\ell +1} \subset C_{\ell}$.

Introduce the quotient group  $Q_{\ell} = \G/C_{\ell}$ with identity element $e_{\ell} \in Q_{\ell}$. There are natural quotient maps $q_{\ell+1} \colon Q_{\ell +1} \to Q_{\ell}$, and we can form the inverse limit group
  \begin{eqnarray} 
\whGamma_{\infty} & \equiv &  \lim_{\longleftarrow} ~ \{ q_{\ell +1} \colon Q_{\ell +1} \to Q_{\ell}  \mid \ell \geq 0 \} \label{eq-invgroup}\\
& = &  \{(g_{\ell}) = (g_0, g_1, \ldots )    \mid g_{\ell} \in Q_{\ell} ~ , ~ q_{\ell +1 }(g_{\ell + 1}) =  g_{\ell} ~ {\rm for ~ all} ~ \ell \geq 0 ~\} ~ \subset \prod_{\ell \geq 0} ~ \G_{\ell} ~ , \label{eq-coordinates}
\end{eqnarray}
which is a Cantor space with the Tychonoff topology. The left actions of $\G$ on the spaces $X_{\ell} = \G/\G_{\ell}$ induce    a minimal  equicontinuous action of $\whGamma_{\infty}$ on $X_{\infty}$, again denoted by  $\whPhi \colon \whGamma_{\infty} \times X_{\infty} \to X_{\infty}$. Note that the isotropy group of the   action of $Q_{\ell} = \G_{\ell}/C_{\ell}$ at the identity coset  in $X_{\ell}=  \G/\G_{\ell}$ is the subgroup $D_{\ell} = \Gamma_{\ell}/C_{\ell}$.

Denote the points in $\whGamma_{\infty}$ by 
$\whg = (g_{\ell}) \in \whGamma_{\infty}$ where $g_{\ell} \in Q_{\ell}$. There is a natural basepoint $\whe_{\infty} \in \whGamma_{\infty}$ given by the cosets of the identity element $e \in \G$, so $\whe_{\infty} = (e_{\ell})$ where $e_{\ell} = e C_{\ell} \in Q_{\ell}$ is the identity element in $Q_{\ell}$. 

For each $\ell \geq 0$, let $\Pi_{\ell} \colon \whGamma_{\infty} \to Q_{\ell}$ denote the projection onto the $\ell$-th factor in \eqref{eq-invgroup}, so in the coordinates of \eqref{eq-coordinates}, we have $\Pi_{\ell}(\whg) = g_{\ell} \in Q_{\ell}$. 
The maps $\Pi_{\ell}$ are continuous for the profinite topology on $\whGamma_{\infty}$, so the pre-images of points in $Q_{\ell}$ are clopen subsets. In particular, the  fiber of $\Pi_{\ell} \colon \whGamma_{\infty} \to Q_{\ell}$  over   $e_{\ell}$ is the normal subgroup
\begin{equation}\label{eq-opennbhds}
\whC_{\ell} = \Pi_{\ell}^{-1}(e_{\ell}) =  \{(g_{i})  \in \whGamma_{\infty}  \mid  g_{i} \in C_{i} ~ , ~ 0 \leq i \leq \ell \} \ . 
\end{equation}

The collection $\{\whC_{\ell} \mid \ell \geq 1\}$ forms a basis of   clopen neighborhoods of $\whe_{\infty} \in \whGamma_{\infty}$. That is, for each clopen set $\whU \subset \whGamma_{\infty}$ with $\whe_{\infty} \in \whU$, there exists $\ell_0 > 0$ such that $\whC_{\ell} \subset \whU$ for all $\ell \geq \ell_0$.
\begin{thm}\cite[Theorem~4.4]{DHL2016a}\label{thm-fundamentaliso}
There is an   isomorphism $\whtau \colon \fG(\Phi) \to \whGamma_{\infty}$ which conjugates the profinite action
$(\fX, \fG(\Phi), \whPhi)$ with the profinite action
$(X_{\infty}, \whGamma_{\infty}, \whPhi)$. In particular, $\whtau$ 
identifies the isotropy group $\fD(\Phi)$ with the inverse limit subgroup
\begin{equation}\label{eq-discformula}
D_{\infty}  = \varprojlim \ \{q_{\ell +1} \colon \G_{\ell +1}/C_{\ell+1} \to \G_{\ell}/C_{\ell} \mid \ell \geq 0\} \subset \whGamma_{\infty}~ .
\end{equation}
\end{thm}

The maps $q_{\ell +1}$ in the formula \eqref{eq-discformula} need not be surjections, and thus the calculation of the inverse limit $D_{\infty}$   can involve some subtleties. For example, it is possible that each group $Q_{\ell}$ is non-trivial for $\ell > 0$, and yet $D_{\infty}$ is the trivial group.

\subsection{Equivalence of Cantor actions}\label{subsec-equivalence}

We next recall the   notions of equivalence of  Cantor actions. 
The first and strongest   is that  of 
  \emph{isomorphism}, which   is a   generalization  of the   notion of conjugacy of topological actions. For $\G = \mZ$, isomorphism corresponds to the notion of ``flip conjugacy'' introduced    in the work of Boyle and Tomiyama \cite{BoyleTomiyama1998}.
The definition below also appears in the papers     \cite{CortezMedynets2016,HL2019,Li2018}.

 \begin{defn} \label{def-isomorphism}
Cantor actions $(\fX_1, \G_1, \Phi_1)$ and $(\fX_2, \G_2, \Phi_2)$ are said to be \emph{isomorphic}  if there is a homeomorphism $h \colon \fX_1 \to \fX_2$ and a group isomorphism $\Theta \colon \G_1 \to \G_2$ so that 
\begin{equation}\label{eq-isomorphism}
\Phi_1(g) = h^{-1} \circ \Phi_2(\Theta(g)) \circ h   \in   \Homeo(\fX_1) \   \textrm{for  all} \ g \in \G_1 \ .
\end{equation}
 \end{defn}

The notion of \emph{return equivalence} for Cantor actions is  weaker than   isomorphism, and is natural when considering the dynamical properties of Cantor systems which should be independent of the restriction of the action to a clopen cross-section.

For a minimal equicontinuous Cantor action $(\fX, \G, \Phi)$ and   an adapted set $U \subset \fX$, by a small abuse of  notation, we use $\Phi_U$ to denote both the restricted action $\Phi_U \colon \G_U \times U \to U$ and the induced quotient action $\Phi_U \colon \cH_U \times U \to U$ for $\cH_U = \Phi(G_U) \subset \Homeo(U)$. Then $(U, \cH_U, \Phi_U)$ is called the \emph{holonomy action} for $\Phi$.
 
   \begin{defn}\label{def-return}
Two minimal equicontinuous Cantor  actions $(\fX_1, \G_1, \Phi_1)$ and $(\fX_2, \G_2, \Phi_2)$ are  \emph{return equivalent} if there exists 
  an adapted set $U_1 \subset \fX_1$ for the action $\Phi_1$   and  
  an adapted set $U_2 \subset \fX_2$ for the action $\Phi_2$,
such that   the  restricted actions $(U_1, \cH_{1,U_1}, \Phi_{1,U_1})$ and $(U_2, \cH_{2,U_2}, \Phi_{2,U_2})$ are isomorphic.
\end{defn}
If the actions $\Phi_1$ and $\Phi_2$ are isomorphic in the sense of Definition~\ref{def-isomorphism}, then they are return equivalent with   $U_1 = \fX_1$ and $U_2 = \fX_2$. However, the notion of return equivalence is weaker even for this case, as the conjugacy is between the holonomy groups $\cH_{1,\fX_1}$ and $\cH_{2,\fX_2}$, and not the groups $\G_1$ and $\G_2$.

\subsection{Locally quasi-analytic}\label{subsec-lqa}

 The quasi-analytic property for Cantor actions  was introduced by
 {\'A}lvarez L{\'o}pez and  Candel  in  \cite[Definition~9.4]{ALC2009} as a generalization of the notion of a \emph{quasi-analytic action} studied by Haefliger for actions of pseudogroups of real-analytic diffeomorphisms.  The authors introduced a local form of the quasi-analytic property  in \cite{DHL2016c,HL2019}: 
  
\begin{defn} \cite[Definition~2.1]{HL2019} \label{def-LQA} A topological action       $(\fX,\G,\Phi)$ on a metric Cantor space $\fX$,  is   \emph{locally quasi-analytic}  if there exists $\e > 0$ such that for any non-empty open set $U \subset \fX$ with $\diam (U) < \e$,  and  for any non-empty open subset $V \subset U $, and elements $g_1 , g_2 \in \G$
 \begin{equation}\label{eq-lqa}
  \text{if the restrictions} ~~ \Phi(g_1)|V = \Phi(g_2)|V, ~ \text{ then}~~ \Phi(g_1)|U = \Phi(g_2)|U. 
\end{equation}
 The action is said to be \emph{quasi-analytic} if \eqref{eq-lqa} holds for $U=\fX$.
\end{defn}

In other words, $(\fX,\G,\Phi)$ is locally quasi-analytic if for every $g \in \G$, the homeomorphism $\Phi(g)$ has unique extensions on the sets of diameter $\e>0$ in $\fX$, with $\e$ uniform over $\fX$. We note that an effective action $(\fX,\G,\Phi)$, for a countable group $\G$,  is topologically free if and only if it is quasi-analytic. 

Recall that a group $\G$ is \emph{Noetherian}  \cite{Baer1956} if every increasing chain of subgroups has a maximal element. Equivalently, a group is Noetherian if every subgroup of $\G$ is finitely generated.

\begin{thm}\label{thm-noetherian} \cite[Theorem~1.6]{HL2020}
Let    $\G$ be a Noetherian group. Then   a  minimal equicontinuous Cantor action $(\fX,\G,\Phi)$   is locally quasi-analytic.
\end{thm}
 A finitely-generated nilpotent group is Noetherian, so as a corollary we obtain that all Cantor actions by finitely-generated nilpotent groups are locally quasi-analytic.   
 
 The notion of an LQA Cantor action extends to the case of a profinite group action $\whPhi \colon \fG \times \fX \to \fX$. 
     
       \begin{defn} \label{def-stable2}
       Let  $(\fX,\G,\Phi)$ be a Cantor action, and  $\whPhi \colon \fG \times \fX \to \fX$ the induced profinite action.   We say that the action is \emph{stable} if the induced profinite action $(\fX, \fG(\Phi), \whPhi)$  is locally quasi-analytic, and   is said to be \emph{wild} otherwise.
 \end{defn}
A profinite completion  $\fG$ of a Noetherian group $\G$ need not be  Noetherian, as can be seen for the example of $\G = \mZ$, and    $\whGamma$  the full profinite completion of $\mZ$.  In particular, the assumption that the action $(\fX,\G,\Phi)$ is locally quasi-analytic does not imply that the action is stable.

 \subsection{Type and typeset for  Cantor actions}

A  Steinitz number $\xi$ can be written uniquely as the formal product over the set of primes $\Pi$, 
\begin{equation}\label{eq-defsteinitz}
\xi  = \prod_{p \in \Pi} ~ p^{\chi_{\xi}(p)} 
\end{equation}
where the \emph{characteristic function} $\chi_{\xi} \colon \Pi \to \{0,1,\ldots, \infty\}$ counts the multiplicity with which a prime $p$ appears in the infinite product $\xi$.

\begin{defn}\label{defn-type}
Two Steinitz numbers  $\xi$ and $\xi'$ are said to be \emph{asymptotically equivalent} if there exists finite integers $m, m' \geq 1$ such that $m \cdot \xi = m' \cdot \xi'$, and we then write $\xi \mor \xi'$.
 
A \emph{type} is an asymptotic equivalence class of Steinitz numbers. The type associated to a Steinitz number  $\xi$ is   denoted by $\tau[\xi]$.
\end{defn}

  In terms of their    characteristic functions $\chi_1, \chi_2$, we have $\xi \mor \xi'$ if and only if the following conditions are satisfied:
\begin{itemize}
\item $\chi_1(p) = \chi_2(p)$ for all but finitely many primes $p \in \Pi$,
\item $\chi_1(p) = \infty$ if and only iff $\chi_1(p) = \infty$ for all primes $p \in \Pi$.
\end{itemize}

Given two types, $\tau$ and $\tau'$ we write $\tau \leq \tau'$, if there exists representatives $\xi \in \tau$ and $\xi' \in \tau'$ such that their characteristic functions satisfy $\chi_{\xi}(p) \leq \chi_{\xi'}(p)$ for all primes $p \in \Pi$.

  \begin{defn}\label{def-primespectrum}
Let  $\pi$ denote the set of   primes. Given  $\xi = \prod_{p \in \pi} \ p^{\chi_{\xi}(p)}$, define:
\begin{eqnarray*}
\pi(\xi) ~ & = & ~  \{ p \in \pi \mid \chi_{\xi}(p) >0   \}   \ , ~ \emph{the prime spectrum of} \ \xi ,  \label{eq-primespectrum}\\
\pi_f(\xi) ~ & = & ~  \{ p \in \pi \mid 0 < \chi_{\xi}(p) < \infty \} \ , ~  \emph{the finite prime spectrum of} \ \xi ,  \label{eq-finiteprimespectrum} \\
\pi_{\infty}(\xi) ~ & = & ~  \{ p \in \pi \mid  \chi_{\xi}(p) = \infty \} \ , ~  \emph{the infinite prime spectrum of} \ \xi \ .\label{eq-infiniteprimespectrum}
\end{eqnarray*}
   \end{defn}
  Note that if $\xi \mor \xi'$, then $\pi_{\infty}(\xi) = \pi_{\infty}(\xi')$. The property that $\pi_f(\xi)$ is an \emph{infinite} set is also preserved by asymptotic equivalence of Steinitz numbers.

 We define the type of a Cantor action $(X_{\infty}, \G, \Phi)$ defined by  
 a   chain of finite index subgroups, $\cG = \{\G = \G_0 \supset \G_1  \supset \cdots\}$. Let $C_{\ell} \subset \G_{\ell}$ denote the normal core of $\G_\ell$ in $\G$.

\begin{defn}\label{def-typegroup}
Let $(X_{\infty}, \G, \Phi)$ be a minimal equicontinuous Cantor action defined by a group chain $\cG$. The type $\tau[X_{\infty}, \G, \Phi]$ of the action is the equivalence class of the Steinitz order
\begin{equation}\label{eq-actionsteinitzorder}
\xi(X_{\infty}, \G, \Phi) = \lcm \{ \# X_{\ell} = \#(\G/\G_{\ell}) \mid \ell > 0 \}  \ .
\end{equation}
\end{defn}

Finally, we note the following result:
\begin{thm}  \cite[Theorem~1.9]{HL2023a} \label{thm-welldefined}
Let $(\fX, \G, \Phi)$ be a Cantor action. The Steinitz order $\xi(\fX, \G, \Phi)$ is defined to be the Steinitz order for  an algebraic model $(X_{\infty}, \G, \Phi)$ of the action, which does not depend upon the choice of an algebraic model. Moreover,  the type $\tau[\fX, \G, \Phi]$ depends only on the return equivalence class of the action.
\end{thm}

\subsection{Type for profinite groups}\label{subsec-protype}

 The \emph{Steinitz order} $\Pi[\fG]$ of a   profinite group $\fG$ is  defined by the supernatural number associated to a presentation of $\fG$ as an inverse limit of finite groups (see    \cite[Chapter 2]{Wilson1998} or \cite[Chapter~2.3]{RZ2000}).
 The Steinitz order appears in the study of   analytic representations of profinite groups associated to groups acting on rooted  trees, see for example   \cite{Kionke2019}.  
 
Recall that for  a profinite group $\fG$, an open subgroup $\fU \subset \fG$ has finite index  \cite[Lemma 2.1.2]{RZ2000}. 
\begin{defn}\label{def-steinitzorderaction}
Let $(\fX, \G, \Phi)$ be a minimal equicontinuous Cantor action, with choice of a basepoint $x \in \fX$. The \emph{Steinitz orders} of the action  are defined as follows:
 \begin{enumerate}
\item $\xi(\fG(\Phi)) =  \lcm \{\# \ \fG(\Phi)/\fN  \mid \fN \subset \fG(\Phi)~ \text{open normal subgroup}\}$, 
\item $\xi(\fD(\Phi)) =  \lcm \{\# \ \fD(\Phi)/(\fN \cap \fD(\Phi)) \mid \fN \subset \fG(\Phi)~ \text{open normal subgroup}\}$, 
\item $\xi(\fG(\Phi) : \fD(\Phi)) =  \lcm \{\# \ \fG(\Phi)/(\fN \cdot \fD(\Phi))  \mid  \fN \subset \fG(\Phi)~ \text{open normal subgroup}\}$.
\end{enumerate}
\end{defn}
The Steinitz orders  satisfy the Lagrange identity, where the multiplication is taken in the sense of supernatural numbers, 
 \begin{equation}\label{eq-productorders}
\xi(\fG(\Phi))  = \xi(\fG(\Phi) : \fD(\Phi)) \cdot \xi(\fD(\Phi))  \ .
\end{equation}
 and thus we always have $\tau[\fD(\Phi)] \leq \tau[\fG(\Phi)]$. The following is a direct consequence of the definitions:
 \begin{thm}
 Let $(\fX, \G, \Phi)$ be a Cantor action. Then there is equality of Steinitz orders, $\xi(\fX, \G, \Phi) = \xi(\fG(\Phi) : \fD(\Phi))$.
 \end{thm}

 \section{Nilpotent actions}\label{sec-nilpotent}

In this section, we apply the notion of  the Steinitz order of a nilpotent Cantor action to the study of  its dynamical properties. 
The proof of Theorem~\ref{thm-main1}    is based on  the special properties of the profinite completions of nilpotent groups, 
in particular the uniqueness of their Sylow $p$-subgroups, and on the relation of this algebraic property with the dynamics of the action.

\subsection{Noetherian groups}\label{subsec-noetherian}
 A countable group $\G$ is said to be \emph{Noetherian} \cite{Baer1956} if every increasing chain of   subgroups $\{H_i \mid i \geq 1 \}$ of $\G$ has a maximal element $H_{i_0}$. The group $\mZ$ is Noetherian;  a finite product of Noetherian groups is Noetherian; and a subgroup and quotient group of a Noetherian group is Noetherian.  Thus, a finitely-generated nilpotent group is Noetherian. 
  
The notion of a Noetherian group has a generalization which is useful for the study of actions of  profinite groups.  
\begin{defn}  \cite[page 153]{Wilson1998} \label{def-noetherian} 
A profinite group $\fG$ is said to be \emph{topologically Noetherian} if every increasing chain of \emph{closed} subgroups $\{\fH_i \mid i \geq 1 \}$ of $\fG$ has a maximal element $\fH_{i_0}$.
\end{defn}
   
   We illustrate this concept with two canonical examples of profinite completions of $\mZ$. 
   
    \begin{ex}\label{ex-stableZ}
 {\rm 
 Let   $\whmZ_p$ denote the $p$-adic integers, for $p$ a prime. That is, $\whmZ_p$ is the completion of $\mZ$ with respect to the chain of subgroups 
 $\cG = \{\G_{\ell} = p^{\ell} \mZ \mid \ell \geq 1\}$. The closed subgroups of  $\whmZ_p$ are given by   $p^i \cdot \whmZ_p$ for some fixed $i > 0$, hence satisfy the ascending chain property  in Definition \ref{def-noetherian}. 
 }
\end{ex}

 \begin{ex}\label{ex-wildZ}
 {\rm 
 Let $\whpi = \{p_i \mid i \geq 1\}$ be an infinite collection of distinct primes. Define an increasing chain of subgroups of $\mZ$, where   $\cG_{\whpi} = \{\G_{\ell} = p_1p_2 \cdots p_{\ell} \mZ \mid \ell \geq 1\}$. Let $\whmZ_{\whpi}$ be the completion of $\mZ$ with respect to the chain $\cG_{\whpi}$. Then we have a topological isomorphism
 \begin{equation}
\whmZ_{\whpi} \cong \prod_{i \geq 1} \ \mZ/p_i \mZ \ .
\end{equation}
Let $H_{\ell} = \mZ/p_1\mZ \oplus \cdots \oplus \mZ/p_{\ell} \mZ$ be the direct sum of the first $\ell$-factors. Then $\{H_{\ell} \mid \ell \geq 1\}$ is an   increasing chain of   subgroups of $\whmZ_{\whpi}$   which does not stabilize, so $\whmZ_{\whpi}$ is not topologically Noetherian.
}
\end{ex}
  These two examples illustrate  the   idea behind the proof of the following result. 
\begin{prop}\label{prop-nilpNoetherian}
Let $\G$ be a finitely generated nilpotent group, and let $\whGamma$ be a profinite completion of $\G$.
Then $\whGamma$ is topologically Noetherian if and only if   the prime spectrum $\pi(\xi(\whGamma))$ is finite.
\end{prop}
 \proof 
First, recall some basic facts about profinite groups. (See for example,  \cite[Chapter~2]{Wilson1998}.) For a prime $p$,  a finite group $H$ is a $p$-group if every element of $H$ has order a power of $p$. A profinite group $\fH$ is a pro-$p$-group if $\fH$ is the inverse limit of finite $p$-groups. A Sylow $p$-subgroup $\fH \subset \fG$ is a maximal pro-$p$-subgroup \cite[Definition~2.2.1]{Wilson1998}. 

A profinite group $\fG$ is \emph{pro-nilpotent} if it is the inverse limit of finite nilpotent groups.  For example, if $\fG$ is a profinite completion of a nilpotent group $\G$, then $\fG$ is pro-nilpotent. 

The group $\fG$ is topologically finitely generated if it contains a dense subgroup $\G \subset \fG$ where $\G$ is finitely generated. The completion $\fG(\Phi)$ associated to a Cantor action $(\fX, \G, \Phi)$ with $\G$ finitely generated is topologically finitely generated.

Assume that  $\fG$ is pro-nilpotent, then for each prime $p$, there is a unique Sylow $p$-subgroup of $\fG$, which is normal in $\fG$ (see \cite[Proposition~2.4.3]{Wilson1998}). Denote this group by $\fG_{(p)}$. Moreover, $\fG_{(p)}$ is non-trivial if and only if 
$p \in \pi(\xi(\fG))$. We use the   following result for pro-nilpotent groups, which is a consequence of \cite[Proposition~2.4.3]{Wilson1998}. 
\begin{prop}\label{prop-factorization}
Let $\fG$ be a profinite completion of a finitely-generated nilpotent group $\G$. Then there is a topological isomorphism
\begin{equation}\label{eq-primeSylowdecomp}
\fG \cong \prod_{p \in \pi(\xi(\fG))} \ \fG_{(p)} \ .
\end{equation}
\end{prop}

 From the isomorphism \eqref{eq-primeSylowdecomp} it follows immediately that if the prime spectrum $\pi(\xi(\fG))$ is infinite, then $\fG$ is not topologically Noetherian. To see this, list $\pi(\xi(\fG)) = \{p_i \mid i = 1,2, \ldots \}$, then we obtain an infinite strictly  increasing chain of closed subgroups,
 $$\fH_{\ell} = \prod_{i=1}^{\ell} \ \fG_{(p_i)} \ . $$
 If the prime spectrum $\pi(\xi(\fG))$ is finite,  then the isomorphism   \eqref{eq-primeSylowdecomp} reduces the proof that $\fG$ is topologically Noetherian to the case of showing that if $\fG$ is topologically finitely generated, then each of its Sylow $p$-subgroups is Noetherian. The group  $\fG_{(p)}$ is nilpotent and topologically finitely generated, so we can use  the lower central series for $\fG_{(p)}$ and induction to reduce to the case where 
 $\fH$ is a topologically finitely-generated abelian pro-$p$-group, and so is isomorphic to a finite product of $p$-completions of $\mZ$, which are topologically Noetherian.
  
 Observe  that   a   profinite completion $\fG$ of a finitely generated nilpotent group $\G$ is a topologically finitely-generated nilpotent group, and we apply the above remarks.
\endproof

\begin{cor}\label{cor-vnilpotentNoetherian}
Let $\G$ be a virtually nilpotent group; that, is there exists a finitely-generated nilpotent subgroup $\G_0 \subset \G$ of finite index.
Then   a    profinite completion $\fG$ of $\G$  is topologically Noetherian if and only if   its prime spectrum $\pi(\xi(\fG))$ is finite.
\end{cor}
\proof
We can assume that $\G_0$ is a normal subgroup of $\G$, then its closure $\fG_0 \subset \fG$ satisfies the hypotheses of Proposition~\ref{prop-nilpNoetherian}, and the Steinitz orders satisfy $\xi(\fG_0) \mor \xi(\fG)$.
As $\fG_0$ is topologically Noetherian if and only if $\fG$ is topologically Noetherian, the claim follows.
\endproof

\subsection{Dynamics of Noetherian groups}\label{subsec-Noetheriandynamics}
We   relate the topologically Noetherian property of a profinite group with the dynamics of a Cantor action of the group, to obtain the proof of Theorem~\ref{thm-main1}.  We first give  the profinite analog of     \cite[Theorem~1.6]{HL2020}. We follow the  outline of its proof in  \cite{HL2020}.
  
\begin{prop}\label{prop-NLQA}
Let $\fG$ be a topologically Noetherian group.
Then     a minimal equicontinuous action   $(\fX,\fG,\whPhi)$  on a Cantor space $\fX$   is locally quasi-analytic.
\end{prop}
\proof
The closure $\fG(\Phi)  \subset \Homeo(\fX)$,   so     the action $\whPhi$ of $\fG(\Phi)$ is effective.  Suppose that the action $\whPhi$ is not locally quasi-analytic, then there exists an infinite properly decreasing chain of clopen subsets of $\fX$,
$\{U_1 \supset U_2 \supset \cdots \}$, which satisfy the following properties, for all $\ell \geq 1$:
\begin{itemize}
\item $U_{\ell}$ is adapted to the action $\whPhi$ with isotropy subgroup $\fG_{U_{\ell}} \subset \fG$;
\item there is a closed subgroup $K_{\ell} \subset \fG_{U_{\ell+1}}$ whose restricted action to $U_{\ell +1}$ is trivial, but the restricted action of $K_{\ell}$ to $U_{\ell}$ is effective.
\end{itemize}
It follows that we obtain a properly increasing chain of closed subgroups $\{K_1 \subset K_2 \subset \cdots\}$ in $\fG$, which contradicts the assumption that $\fG$ is topologically Noetherian.
\endproof

\proof[Proof   of Theorem~\ref{thm-main1}]
Let $(\fX,\G,\Phi)$ be a  nilpotent Cantor action, and we are given that the prime spectrum $\pi(\xi(\fG(\Phi)))$ is   finite.  
 Then there exists a finitely-generated  nilpotent   subgroup $\G_0 \subset \G$ of finite index, and we can assume without loss  of generality   that $\G_0$ is normal. Let $\fG(\Phi)_0$ be the closure of $\G_0$ in $\fG(\Phi)$.
The group $\fG(\Phi)$ has finite prime spectrum implies that the group $\fG(\Phi)_0$ has finite prime spectrum, and thus by Proposition~\ref{prop-nilpNoetherian} the group $\fG(\Phi)_0$ is topologically Noetherian. Let $x \in \fX$, then  it suffices to show that the action of $\G_0$ on the orbit $\fX_0 = \fG(\Phi)_0 \cdot x$ is stable. This reduces the proof to showing the claim when  $\G$ is   nilpotent. Then  the profinite closure $\fG(\Phi)$   is also nilpotent, and we have a profinite action $(\fX, \fG(\Phi) ,\whPhi)$.

 Suppose that the action $\whPhi$ is not locally quasi-analytic, then   there exists an increasing  chain of closed subgroups $K_{\ell} \subset \fD(\Phi)$  where $K_{\ell}$ acts trivially on the clopen subset $U_{\ell} \subset \fX$.   As $\fD(\Phi)$ is a closed subgroup of $\fG(\Phi)$, the increasing chain $\{K_{\ell} \mid \ell > 0\}$ consists of closed subgroups of $\fG(\Phi)$.  This contradicts the fact that $\fG(\Phi)$ is topologically Noetherian.   Hence, the action $\whPhi$ must be locally quasi-analytic. That is,  the action $(\fX,\G,\Phi)$ is stable.
  \endproof

 \section{Basic examples}\label{sec-examples}

In this section, we construct two basic  examples of nilpotent Cantor actions.  These   examples illustrate the principles behind the subsequent more complex constructions in Section \ref{sec-prescribed}, which are used to prove Theorems~\ref{thm-main3} and \ref{thm-main2}.

 The integer Heisenberg group  is the simplest non-abelian nilpotent group, and it can be  represented as the upper triangular matrices in ${\rm GL(3,\mZ)}$. That is,
\begin{equation}\label{eq-cH}
\G =   \left\{  \left[ {\begin{array}{ccc}
   1 & a & c\\
   0 & 1 & b\\
  0 & 0 & 1\\
  \end{array} } \right] \mid a,b,c  \in \mZ\right\} .
\end{equation}
 We denote a $3 \times 3$ matrix in $\G$ by the coordinates as $(a,b,c)$.

 \begin{ex}\label{ex-trivial}
 {\rm
 A \emph{renormalizable Cantor action}, as defined in \cite{HLvL2020}, can be constructed   from the group chain defined by a proper self-embedding of a non-abelian group $\G$ into itself.

For a  prime $p \geq 2$, define the self-embedding $\vp_p \colon \G \to \G$ by  
$\vp(a,b,c) = (pa, pb, p^2c)$. Then define a group chain in $\G$ by setting 
$$\G_{\ell} = \vp_p^{\ell}(\G) = \{(p^{\ell} a, p^{\ell}b, p^{2\ell}c) \mid a,b,c \in \mZ\} \quad, \quad \bigcap_{\ell > 0} \ \G_{\ell} = \{e\} \ .$$
For $\ell > 0$, the normal core for $\G_{\ell}$ is given by
$C_{\ell} = {\rm core}(\G_{\ell})  = \{(p^{2\ell} a, p^{2\ell} b, p^{2\ell} c) \mid a,b,c \in \mZ\}$, and so the quotient group  
$Q_{\ell} = \G/C_{\ell} \cong \{( \oa, \ob, \oc) \mid \oa, \ob, \oc\in \mZ/p^{2\ell}\mZ \}$.
The profinite group $\whGamma_{\infty}$ is the inverse limit of the quotient groups $Q_{\ell}$ so we have
$\ds \whGamma_{\infty} =   \{(\wha, \whb, \whc) \mid \wha, \whb, \whc\in \widehat{\mZ}_{p^2} \}$.
   Thus,   $\xi(\whGamma) = \{p^{\infty}\}$. Even though the quotient groups $\G_{\ell}/C_{\ell}$ are all non-trivial, for this action the inverse limit $\fD(\Phi)$ is the trivial group. This follows from the fact that there are inclusions 
   $$\G_{2\ell} = \{ (p^{2\ell}a, p^{2\ell}b,p^{4\ell}c) \mid a,b,c, \in \mZ\} \subset C_\ell =  \{(p^{2\ell} a, p^{2\ell} b, p^{2\ell} c) \mid a,b,c \in \mZ\} \ . $$
The triviality of  $\fD(\Phi)$ implies that there is an equivalent  group chain for the action \cite{DHL2016a} which can   be chosen so that every subgroup in the chain is normal in $\G$. 
 }
 \end{ex}
 
  \begin{ex}\label{ex-2primes}
 {\rm
For distinct  primes $p, q \geq 2$, define the self-embedding $\vp_{p,q} \colon \G \to \G$ by  
$\vp(a,b,c) = (pa, qb, pqc)$. Then define a group chain in $\G$ by setting 
$$\G_{\ell} = \vp_{p,q}^{\ell}(\G) = \{(p^{\ell} a, q^{\ell}b, (pq)^{\ell}c) \mid a,b,c \in \mZ\} \quad, \quad \bigcap_{\ell > 0} \ \G_{\ell} = \{e\} \ .$$

For $\ell > 0$, the normal core for $\G_{\ell}$ is given by
$C_{\ell} = {\rm core}(\G_{\ell})  = \{((pq)^{\ell} a, (pq)^{\ell} b, (pq)^{\ell} c) \mid a,b,c \in \mZ\}$, and so  the quotient group  
$Q_{\ell} = \G/C_{\ell} \cong \{( \oa, \ob, \oc) \mid \oa, \ob, \oc\in \mZ/(pq)^{\ell}\mZ \}$.
The profinite group $\whGamma_{\infty}$ is the inverse limit of the quotient groups $Q_{\ell}$ so we have
$\ds \whGamma_{\infty} =   \{(\wha, \whb, \whc) \mid \wha, \whb, \whc\in \widehat{\mZ}_{pq} \}$.
Thus, $\xi(\whGamma_{\infty}) = \{ p^{\infty}, q^{\infty}\}$, and $D_{\infty}$ is the inverse limit of the finite groups $\G_{\ell}/C_{\ell}$ by \eqref{eq-discformula},  so   $D_{\infty} \cong \widehat{\mZ}_{q} \times \widehat{\mZ}_{p}$.
 }
 \end{ex}

 \section{Nilpotent actions with prescribed spectrum}\label{sec-prescribed}
 
 In this section, we construct stable actions of the discrete Heisenberg group with prescribed prime spectrum, proving Theorem \ref{thm-main3}. Then we construct examples of wild nilpotent Cantor actions, proving Theorem \ref{thm-main2} from which we deduce  Corollary \ref{cor-main4}.
 For simplicity, our examples all use the Heisenberg group represented by $3 \times 3$ matrices. Of course, these examples can be generalized to the integer upper triangular matrices in all dimensions, where there is much more freedom in the choices made in the construction. The   calculations become correspondingly more tedious, but yield analogous results. It seems reasonable to expect that similar constructions can be made for any finitely-generated torsion-free nilpotent (non-abelian) group $\G$.   That is,   there are always group chains in $\G$ which yield wild nilpotent Cantor   actions.  
 
Let  $\G \subset {\rm GL}(3,\mathbb{Z})$ denote  the discrete Heisenberg group, given by formula \eqref{eq-cH}. The basis for the constructions below is the structure theory for nilpotent group completions  in  Proposition~\ref{prop-factorization}, in particular the formula \eqref{eq-primeSylowdecomp}. Given sets of primes $\pi_f$ and $\pi_{\infty}$, we embed an   infinite product  of finite actions, as in Example~\ref{ex-toy}, into   a profinite completion $\whGamma_{\infty}$ of $\G$, and thus obtain  a nilpotent Cantor action $(X_{\infty}, \G, \Phi)$  on the quotient   space $X_{\infty} = \whGamma_{\infty}/D_{\infty}$.  
 
 \subsection{Basic components of the construction.}\label{ex-toy}

   Fix a prime $p \geq 2$. 
   
   For $n \geq 1$ and  $0 \leq k < n$, we have the following finite groups: 
 \begin{equation}
G_{p,n} =   \left\{  \left[ {\begin{array}{ccc}
   1 & \oa & \oc\\
   0 & 1 & \ob\\
  0 & 0 & 1\\
  \end{array} } \right] \mid \oa,\ob,\oc  \in \mZ/p^{n}\mZ\right\} ~ , ~
  H_{p,n,k} =   \left\{  \left[ {\begin{array}{ccc}
   1 & p^k \oa & 0\\
   0 & 1 & 0\\
  0 & 0 & 1\\
  \end{array} } \right] \mid \oa  \in \mZ/p^{n}\mZ\right\}
\end{equation}
Note that $\#[G_{p,n}] = p^{3n}$ and $\#[H_{p,n,k}] = p^{n-k}$.

Let   $\ox = (1,0,0), \oy = (0,1,0), \oz= (0,0,1) \in G_{p,n}$, then  
$\ox \cdot \oy \cdot \ox^{-1} = \oy \oz$ and $\ox \cdot \oz \cdot \ox^{-1} = \oz$.  That is, the adjoint action of $\ox$ on the ``plane'' in the  $(\oy,\oz)$-coordinates   is a ``shear'' action along the $\oz$-axis, and the adjoint action of $\ox$ on the $\oz$-axis fixes   all points on the $\oz$-axis.

Set $X_{p,n,k} = G_{p,n}/H_{p,n,k}$,  then the isotropy group of the action of $G_{p,n}$ on $X_{p,n,k}$ at the coset  $H_{p,n,k}$  of the identity element is $H_{p,n,k}$. The core subgroup $C_{p,n,k} \subset H_{p,n,k}$ contains elements in $H_{p,n,k}$ which fix every point in $X_{p,n,k}$. The action of $\ox \in H_{p,n,k}$ on the coset space $X_{p,n,k}$ satisfies 
 $$\Phi(\ox)(\oy\, H_{p,n,k}) = \oy \oz\, H_{p,n,k},$$ so the identity is the only element in $G_{p,n}$ which acts trivially on every coset in $X_{p,n,k}$, so $C_{p,n,k}$ is the trivial group. 
 Then $D_{p,n,k} = H_{p,n,k}/C_{p,n,k} = H_{p,n,k}$, and for each $g \in H_{p,n,k}$ its action fixes the cosets of the multiples of $\oz$.

 \subsection{Stable nilpotent actions with finite or infinite prime spectrum}\label{subsec-stableprime}
 
We now prove Theorem \ref{thm-main3} by constructing a family of stable examples with prescribed prime spectra.
 
Let $\pi_f $ and  $\pi_{\infty}$ be two disjoint   collections of primes, with $\pi_f$ a finite set, and $\pi_{\infty}$ a non-empty set.   

Enumerate $\pi_f = \{q_1, q_2, \ldots, q_m\}$, and then choose integers $1 \leq r_i \leq n_i$ for $1 \leq i \leq m$. 

Enumerate $\pi_{\infty} = \{p_1, p_2, \ldots\}$  with the convention (for notational convenience) that if $\ell$ is greater than the number of primes in $\pi_{\infty}$ then we set $p_{\ell} = 1$.
For each $\ell \geq 1$,  define the integers
\begin{eqnarray}
M_{\ell}    & =   &  q_1^{r_1} q_2^{r_2} \cdots q_m^{r_m} \cdot p_1^{\ell} p_2^{\ell} \cdots p_{\ell}^{\ell} \  , \\
N_{\ell}    & =   &  q_1^{n_1} q_2^{n_2} \cdots q_m^{n_m} \cdot p_1^{\ell} p_2^{\ell} \cdots p_{\ell}^{\ell} \  .
\end{eqnarray}
 
For all $\ell \geq 1$, observe that $M_{\ell}$ divides $N_{\ell}$.    

Define a subgroup of the Heisenberg group $\G$, in the coordinates above, 
$$\G_{\ell} = \{ (a M_{\ell},b N_{\ell} ,c N_{\ell}) \mid a,b,c \in \mZ \}.$$
Its  core  subgroup is given by 
$C_{\ell} =  \{ (a N_{\ell},b N_{\ell} ,c N_{\ell}) \mid a,b,c \in \mZ \}$. 
Observe that
$$\mZ/N_{\ell} \mZ \cong \mZ/q_1^{n_1}\mZ \oplus \cdots \oplus  \mZ/q_m^{n_m}\mZ \oplus \mZ/p_1^{\ell}\mZ \oplus \cdots \oplus \mZ/p_{\ell}^{\ell}\mZ \ .$$
 By Proposition~\ref{prop-factorization},  and in the notation of  Example~\ref{ex-toy},  we have for   $k_i = n_i - r_i$ that
\begin{equation}\label{eq-lqafactors}
 \whGamma_{\infty}~ = ~\lim_{\longleftarrow}\{\G/C_\ell \to \G/C_{\ell-1} \mid \ell \geq 1\} ~ \cong ~ \prod_{i=1}^m \ G_{q_i, n_i} ~ \cdot ~ \prod_{j=1}^{\infty} \ \whGamma_{(p_j)}, \end{equation}
 \begin{equation}\label{eq-lqafactors-1}
 \quad D_{\infty} ~ = ~\lim_{\longleftarrow}\{\G_\ell/C_\ell \to \G_{\ell-1}/C_{\ell-1} \mid \ell \geq 1\} ~ \cong ~  \prod_{i=1}^m \ H_{q_i, n_i, k_i} \ . 
\end{equation}
Then the Cantor space  $X_{\infty} = \whGamma_{\infty}/D_{\infty}$ associated to the group chain $\{\G_{\ell} \mid \ell \geq 1\}$ is given by
\begin{equation}\label{eq-lqaspace}
X_{\infty} ~ \cong ~ \prod_{i=1}^m \ X_{q_i, n_i, k_i} ~ \times ~ \prod_{j=1}^{\infty} \ \whGamma_{(p_j)}  \ . 
\end{equation}

In particular,  as the first factor in   \eqref{eq-lqaspace} is a finite product of finite sets, the second factor defines an open neighborhood $$U = \prod_{i=1}^m \ \{x_i\} ~ \times ~ \prod_{j=1}^{\infty} \ \whGamma_{(p_j)}$$ 
where $x_i \in X_{q_i, n_i, k_i}$ is the basepoint given by the coset of the identity element.
That is, $U$ is a clopen neighborhood of the basepoint in $X_{\infty}$.  The isotropy group of $U$ is given by 
\begin{equation}
 \whGamma_{\infty}|U ~ = ~  \prod_{i=1}^m \ H_{q_i, n_i, k_i} ~ \times ~ \prod_{j=1}^{\infty} \ \whGamma_{(p_j)}  \ . 
\end{equation}
The restriction of  $\whGamma_{\infty}|U$ to $U$ is isomorphic to the subgroup  
\begin{equation}
 K|U ~ = ~  \prod_{i=1}^m \ \{\overline{e}_i\} ~ \times ~ \prod_{j=1}^{\infty} \ \whGamma_{(p_j)} ~ \subset ~ \Homeo(U) \ , 
\end{equation}
where $\overline{e}_i \in G_{q_i, n_i}$ is the identity element.  The group $K|U$   acts freely on $U$, and thus the action of $\whGamma_{\infty}$ on $X_{\infty}$ is locally quasi-analytic.
Moreover, the prime spectrum of the action of $\G$ on $X_{\infty}$ is the union $\whpi = \pi_f \cup \pi_{\infty} = \pi(\xi(\whGamma_{\infty}))$.  If $\pi_\infty$ is infinite, then the prime spectrum of the action is infinite.
Note that the group $\G$ embeds into $\whGamma_{\infty}$, since the integers $M_{\ell}$ and $N_{\ell}$ tend to infinity with $\ell$. This completes the proof of Theorem \ref{thm-main3}.

 \subsection{Wild nilpotent actions with infinite prime spectrum}\label{prescribed-wild}

 We now prove Theorem \ref{thm-main2}. We must show that  every infinite set of primes can be realized as the prime spectrum of a wild action of the Heisenberg group $\G$, as defined by \eqref{eq-cH}.
Let $\pi_f $ and  $\pi_{\infty}$ be   disjoint   collections of primes, with $\pi_f$ an infinite set and $\pi_{\infty}$ arbitrary, possibly empty.   

Enumerate $\pi_f = \{q_1, q_2, \ldots\}$ and choose integers $1 \leq r_i < n_i$ for $1 \leq i < \infty$. 

Enumerate $\pi_{\infty} = \{p_1, p_2, \ldots\}$, again with the convention   that if $\ell$ is greater than the number of primes in $\pi_{\infty}$ then we set $p_{\ell} = 1$. 

As in Section \ref{subsec-stableprime}, for  each $\ell \geq 1$,  define the integers
\begin{eqnarray*}
M_{\ell}    & =   &  q_1^{r_1} q_2^{r_2} \cdots q_{\ell}^{r_{\ell}} \cdot p_1^{\ell} p_2^{\ell} \cdots p_{\ell}^{\ell} \  , \\
N_{\ell}    & =   &  q_1^{n_1} q_2^{n_2} \cdots q_{\ell}^{n_{\ell}} \cdot p_1^{\ell} p_2^{\ell} \cdots p_{\ell}^{\ell} \  .
\end{eqnarray*}

For $\ell \geq 1$, define a subgroup of the Heisenberg group $\G$, in the coordinates above, 
\begin{equation}\label{eq-chain2}
\G_{\ell} = \{ (a M_{\ell},b N_{\ell} ,c N_{\ell}) \mid a,b,c \in \mZ \} \ , 
\end{equation}
 Its core  subgroup is given by 
$C_{\ell} =  \{ (a N_{\ell},b N_{\ell} ,c N_{\ell}) \mid a,b,c \in \mZ \}$. For   $k_i = n_i - r_i$ we then have
\begin{equation}\label{eq-lqafactors2}
 \whGamma_{\infty} ~ \cong ~ \prod_{i=1}^{\infty} \ G_{q_i, n_i} ~ \cdot ~ \prod_{j=1}^{\infty} \ \whGamma_{(p_j)} \quad , \quad D_{\infty} ~ \cong ~  \prod_{i=1}^{\infty} \ H_{q_i, n_i, k_i} \ . 
\end{equation}
The Cantor space  $X_{\infty} = \whGamma_{\infty}/D_{\infty}$ associated to the group chain $\{\G_{\ell} \mid \ell \geq 1\}$ is given by
\begin{equation}\label{eq-lqaspace2}
X_{\infty} ~ \cong ~ \prod_{i=1}^{\infty} \ X_{q_i, n_i, k_i} ~ \times ~ \prod_{j=1}^{\infty} \ \whGamma_{(p_j)}  \ . 
\end{equation}
The first factor in   \eqref{eq-lqaspace} is an infinite product of finite sets, so fixing the first $\ell$-coordinates in this product determines a clopen subset of $X_{\infty}$. Let $x_i \in X_{q_i, n_i, k_i}$ denote the  coset of the identity element, which is   the basepoint in $X_{q_i, n_i, k_i}$. Then for  each $\ell \geq 1$, we define a clopen set in $X_{\infty}$ 
\begin{equation}\label{eq-wildclopen}
U_{\ell} = \prod_{i=1}^{\ell} \ \{x_i\} ~ \times ~  \prod_{i=\ell+1}^{\infty} \ X_{q_i, n_i, k_i} ~ \times ~ \prod_{j=1}^{\infty} \ \whGamma_{(p_j)} \ .
\end{equation}
   Recalling the calculations in   Example~\ref{ex-toy}, the subgroup $H_{q_i, n_i, k_i}$ is the isotropy group of the basepoint $x_i \in X_{q_i, n_i, k_i}$. Thus, 
the isotropy subgroup of $U_{\ell}$  for the $\whGamma_{\infty}$-action  is given by the product
\begin{equation}\label{eq-Gammarestr}
 \whGamma_{\infty}|_{U_{\ell}} ~ = ~  \prod_{i=1}^{\ell} \ H_{q_i, n_i, k_i} ~ \times ~  \prod_{i=\ell + 1}^{\infty} \ G_{q_i, n_i} ~ \times ~ \prod_{j=1}^{\infty} \ \whGamma_{(p_j)}  \ . 
\end{equation}
For $j \ne i$, the subgroup $H_{q_i, n_i, k_i}$ acts as the identity on the factors $X_{q_j, n_j, k_j}$ in \eqref{eq-lqaspace2}.
Thus, the  image of  $\whGamma_{\infty}|_{U_{\ell}}$ in $\Homeo(U_{\ell})$  is isomorphic to the subgroup  
\begin{equation}
 Z_{\ell} ~ = ~\whGamma_{\infty}|U_{\ell} ~ = ~  \prod_{i=1}^{\ell} \ \{\overline{e}_i\} ~ \times ~  \prod_{i=\ell + 1}^{\infty} \ G_{q_i, n_i} ~ \times ~ \prod_{j=1}^{\infty} \ \whGamma_{(p_j)}  ~ \subset ~ \Homeo(U_{\ell}) \ , 
\end{equation}
where $\overline{e}_i \in G_{q_i, n_i}$ is the identity element.  

We next show that this action   is not stable; that is, for any $\ell > 0$ there exists a clopen subset $V \subset U_{\ell}$ and non-trivial $\whg \in   Z_{\ell}$ so that the action of $\whGamma_\infty$ restricts to the identity map on $V$. 

We can assume without loss of generality that $V= U_{\ell'}$  for some $\ell' > \ell$. Consider  the restriction map for 
the isotropy subgroup of $Z_{\ell}$ to  $U_{\ell'}$ which is given by 
  $$\rho_{\ell, \ell'} \colon Z_{\ell}|_{U_{\ell'}}  \to Z_{\ell'} \subset \Homeo(U_{\ell'}) \ .$$
We must show that there exists $\ell' > \ell$ such that this map has    a non-trivial kernel.  
 Calculate this map in terms of the product representations above,
  \begin{equation}\label{eq-restriction}
 Z_{\ell}|_{U_{\ell'}}  ~ = ~  \prod_{i=1}^{\ell} \ \{\overline{e}_i\} ~ \times ~  \prod_{i=\ell + 1}^{\ell'} \ H_{q_i, n_i,k_i} ~ \times ~  \prod_{i=\ell' + 1}^{\infty} \ G_{q_i, n_i}~ \times ~ \prod_{j=1}^{\infty} \ \whGamma_{(p_j)}   \ . 
\end{equation}
For $\ell < i \leq \ell'$, the group $H_{q_i, n_i,k_i}$ fixes the point $\prod_{i=1}^{\ell'} \ \{x_i\}$, and acts trivially on 
$\prod_{i=\ell'+1}^{\infty} \ X_{q_i, n_i, k_i}$.
Thus, the kernel of the restriction map contains the second factor in \eqref{eq-restriction}, 
\begin{equation}\label{eq-kernelH}
 \prod_{i=\ell + 1}^{\ell'} \ H_{q_i, n_i,k_i} ~ \subset ~  \ker \left\{ \rho_{\ell, \ell'} \colon Z_{\ell}|_{U_{\ell'}}  \to \Homeo(U_{\ell'})  \right\} \ .
\end{equation}
As this group is non-trivial for all $\ell' > \ell$,  the action of $\whGamma_{\infty}$ on $X_{\infty}$ is not locally quasi-analytic, hence the action of $\G$ on $X_{\infty}$ is wild.
Moreover, the prime spectrum of the action of $\G$ on $X_{\infty}$ equals the union $\whpi = \pi_f \cup \pi_{\infty}$.

We now prove the second part of Theorem~\ref{thm-main2}, showing that choices in the construction above can be made in such a way that the action of $\G$ on a Cantor set is topologically free while the action of $\whGamma_\infty$ is wild, and the prime spectrum is prescribed. 

Choose  an infinite set of distinct primes $\pi_f = \{q_1, q_2, \ldots\}$, and let $\pi_{\infty}$ be empty.

Choose the constants as in Section~\ref{ex-toy},  with $n_i =2$ and $k_i=1$ for all $i \geq 1$. 

Define the Cantor space $X_{\infty}$   by \eqref{eq-lqaspace2}, where the second factor is trivial; that is a point.  The action of $\whGamma_\infty$  is   wild by the calculations in formulas \eqref{eq-Gammarestr} to  \eqref{eq-kernelH}. 

We claim that the action of $\G$ on $X_{\infty}$ is   topologically free. Suppose not, then there exists an open set $U \subset X_{\infty}$ and $g \in \G$ such that the action of $\Phi(g)$ is non-trivial on   $X_{\infty}$ but leaves the set $U$ invariant, and restricts to the identity action on $U$.
The action of $\G$ on $X_{\infty}$ is minimal, so there exists $h \in \G$ with $h \cdot x_{\infty} \in U$. Then $\Phi(h^{-1} g h)(x_{\infty}) = x_{\infty}$ and the action   $\Phi(h^{-1} g h)$ fixes an open neighborhood of $x_{\infty}$. Replacing $g$ with $h^{-1} g h$ we can assume that $\Phi(g)(x_{\infty}) = x_{\infty} \in U$. From the definition \eqref{eq-wildclopen}, the clopen sets
\begin{equation}\label{eq-wildclopen2}
U_{\ell} = \prod_{i=1}^{\ell} \ \{x_i\} ~ \times ~  \prod_{i=\ell+1}^{\infty} \ X_{q_i, 2, 1}  
\end{equation}
form a neighborhood basis at $x_{\infty}$, and thus there exists $\ell > 0$ such that $U_{\ell} \subset U$.

The group $\G$    embeds into $\whGamma_{\infty}$ along the diagonal in the product \eqref{eq-primeSylowdecomp}. That is, we can write    $\ds g = (g,g,\ldots,g) \in    \prod_{i=1}^{\infty} \ G_{q_i, 2}$. The action of $\Phi(g)$ is factorwise, and $\Phi(g)(x_{\infty}) = x_{\infty}$ implies that $\ds g \in D_{\infty} \cong   \prod_{i=1}^{\infty} \ H_{q_i, n_i, k_i}$. The assumption that $\Phi(g)$ fixes the points in $U$ implies that it acts trivially on each factor $X_{q_i, 2, 1}$ for $i > \ell$. As each factor $H_{q_i, 2, 1}$ acts effectively on $X_{q_i, 2, 1}$ this implies that the projection of $g$ to the $i$-th factor group $H_{q_i, 2, 1}$ is the identity for $i > \ell$. This 
  implies that every entry   above the diagonal in the matrix representation of $g$ in \eqref{eq-cH} is divisible by an infinite number of distinct primes $\{q_i \mid i \geq \ell\}$, so by the Prime Factorization Theorem the matrix $g$ is  the identity. 

Alternatively, observe that we have $\ds g \in \prod_{i=1}^{\ell} \ H_{q_i, 2, 1}$. This is a finite product of finite groups, which implies that $g \in \G$ is a   torsion element. However, the Heisenberg group $\G$ is torsion-free, hence $g$ must be the identity.  Thus, the action of $\G$ on $X_{\infty}$ must be topologically free.

Finally, the above construction allows the choice of any infinite subset $\pi_f$  of  distinct primes, and there are an uncountable number such choices which are distinct. 
Thus, by Theorem~1.9 in \cite{HL2023a} there are an uncountable number of topologically-free, wild nilpotent Cantor actions with distinct prime spectrum.
This completes the proof of Theorem~\ref{thm-main2}.

\subsection{Proof of Corollary \ref{cor-main4}}

Consider the family of wild topologically free actions on the Heisenberg group $\G$ with infinite distinct prime spectrum, as  constructed at the end of Section \ref{prescribed-wild}. We show that the uncountable number of infinite choices of $\pi_f$ in this family can be made so that the actions have pairwise disjoint types. 

By Definition \ref{defn-type}, for two Steinitz numbers $\xi$ and $\xi'$ we have that their types are equal, $\tau(\xi) = \tau(\xi')$, if and only if there exist integers $m,m'$ such that $m \cdot \xi = m' \cdot \xi'$. Thus two actions with prime spectra $\pi_f$ and $\pi_f'$ have distinct types if and only if $\pi_f$ and $\pi_f'$ differ by an infinite number of entries. This happens, for instance, if $\pi_f$ and $\pi_f'$ are \emph{almost disjoint infinite sets}, i.e. they are infinite sets with finite intersection.

The set of prime numbers is countable, so the family of infinite almost disjoint subsets of prime numbers is uncountable if and only if the family of infinite almost disjoint subsets of natural numbers is uncountable.  The family of almost disjoint subsets of natural numbers is uncountable by \cite[Corollary 2.3]{Geschke2012}. Since the set of finite subsets of natural numbers is countable, the set of almost disjoint infinite subsets of natural numbers is uncountable.

It follows that the prime spectra of the uncountable family of actions of the Heisenberg group in Theorem \ref{thm-main2} can be chosen so that they form a family of almost disjoint infinite sets. Then their types are pairwise distinct, and by Theorem \ref{thm-welldefined} these actions of the Heisenberg group are pairwise not return equivalent. Therefore, they are pairwise not conjugate.


\end{document}